\documentclass[12pt]{article}
\usepackage{amsmath}
\usepackage{amssymb}
\usepackage{vatola}

\def\q{\quad}
\def\qq{\qquad}
\def\qtq#1{\q\t{#1}\q}
\def\mod#1{\ (\text{\rm mod}\ #1)}
\def\t{\text}
\def\f{\frac}
\def\e{\equiv}
\def\b{\binom}
\def\ap{\langle a\rangle_p}
\def\sls#1#2{(\f{#1}{#2})}
\def\ls#1#2{\big(\f{#1}{#2}\big)}
\def\Ls#1#2{\Big(\f{#1}{#2}\Big)}
\let \pro=\proclaim
\let \endpro=\endproclaim

\begin{document}
\leftline{preprint: May 30, 2016}
\par\q\par\q
 \centerline {\bf
New super congruences involving Bernoulli and Euler polynomials}
\par\q\newline \centerline{Zhi-Hong Sun}
\par\q\newline \centerline{School of
Mathematical Sciences} \centerline{Huaiyin Normal University}
 \centerline{Huaian, Jiangsu 223001, P.R. China}
  \centerline{zhsun@hytc.edu.cn}
\centerline{http://www.hytc.edu.cn/xsjl/szh} \par\q
\newline Let $p>3$ be a prime, and let $a$ be a rational p-adic
integer with $a\not\equiv 0\pmod p$. In this paper we establish
congruences for
$$\sum_{k=1}^{(p-1)/2}\frac{\binom ak\binom{-1-a}k}k,
\quad\sum_{k=0}^{(p-1)/2}k\binom ak\binom{-1-a}k
\quad\text{and}\quad\sum_{k=0}^{(p-1)/2}\frac{\binom
ak\binom{-1-a}k}{2k-1}\pmod {p^2}$$ in terms of Bernoulli and Euler
polynomials. We also give some transformation formulas for
congruences modulo $p^2$.

\par\q
\par\q

\section*{1. Introduction}
\par\q The Bernoulli numbers $\{B_n\}$ and
Bernoulli polynomials $\{B_n(x)\}$ are
 defined by
 $$B_0=1,\ \sum_{k=0}^{n-1}\b nkB_k=0\ (n\ge 2)\qtq{and}
 B_n(x)=\sum_{k=0}^n\b nkB_kx^{n-k}\ (n\ge 0).$$
 The Euler numbers $\{E_n\}$ and Euler polynomials $\{E_n(x)\}$ are
 defined by
$$E_0=1, \ E_n=-\sum_{k=1}^{[n/2]}\b n{2k}E_{n-2k}\ (n\ge 1)
\ \t{and}\ E_n(x)=\f 1{2^n}\sum_{r=0}^n\b nr(2x-1)^{n-r}E_r\ (n\ge
0),$$
 where $[a]$ is the greatest integer not exceeding $a$. In [S6] the author introduced the sequence $\{U_n\}$ given by
$$U_0=1\qtq{and}U_n=-2\sum_{k=1}^{[n/2]}\b n{2k}U_{n-2k}\ (n\ge 1).$$
 It is well
 known that $B_{2n+1}=0$ and $E_{2n-1}=U_{2n-1}=0$ for any positive integer
 $n$. $\{B_n\}$, $\{E_n\}$ and $\{U_n\}$ are important sequences and they have many
interesting properties and applications. See [EMOT], [IR], [MOS] and
[S2-S7].
\par  It is easily seen that ([S8])
$$\aligned&\b{-\f 12}k^2=\f{\b{2k}k^2}{16^k},\ \b{-\f
13}k\b{-\f 23}k=\f{\b{2k}k\b{3k}k}{27^k}, \\&\ \b{-\f 14}k\b{-\f
34}k=\f{\b{2k}k\b{4k}{2k}}{64^k},\ \b{-\f 16}k\b{-\f
56}k=\f{\b{3k}k\b{6k}{3k}}{432^k}.\endaligned\tag 1.1$$  In 2003,
Rodriguez-Villegas [RV] conjectured that for any prime $p>3$,
$$\align &\sum_{k=0}^{p-1}\f{\b{2k}k^2}{16^k}\e
\Ls{-1}p\mod{p^2},\q\sum_{k=0}^{p-1} \f{\b{2k}k\b{3k}k}{27^k}\e\Ls
{-3}p\mod{p^2},\tag 1.2
\\&\sum_{k=0}^{p-1}\f{\b{2k}k\b{4k}{2k}}{64^k}\e \Ls{-2}p\mod{p^2},
\q \sum_{k=0}^{p-1} \f{\b{3k}k\b{6k}{3k}}{432^k}\e
\Ls{-1}p\mod{p^2},\tag 1.3\endalign$$ where $\sls ap$ is the
Legendre symbol. (1.2) and (1.3) were later confirmed by Mortenson
[M1-M2].  Let $\Bbb Z$ be the set of integers. For a prime $p$ let
$\Bbb Z_p$ denote the set of rational numbers whose denominator is
not divisible by $p$. For $a\in\Bbb Z_p$ let
$\ap\in\{0,1,\ldots,p-1\}$ be given by $a\e\ap\mod p$. Let $p>3$ be
a prime, $a\in\Bbb Z_p$, $a\not\e 0\mod p$ and $t=(a-\ap)/p$. In
[S9-S11] the author showed that
$$\align&\q\sum_{k=0}^{p-1}\b ak\b{-1-a}k\e (-1)^{\ap}
+p^2t(t+1)E_{p-3}(-a)\mod{p^3}, \tag 1.4
\\&\q\sum_{k=0}^{p-1}k\b
ak\b{-1-a}k\e (-1)^{\ap}
a(a+1)+p^2t(t+1)(a(a+1)E_{p-3}(-a)-1)\mod{p^3},\tag 1.5
\\&\q\sum_{k=0}^{p-1}\b ak\b{-1-a}k\f 1{2k-1}\tag 1.6
\\&\qq\e -(2a+1)(2t+1)-p^2t(t+1)(4+(2a+1)B_{p-2}(-a))\mod {p^3},
\\&\q\sum_{k=0}^{p-1}\b ak\b{-1-a}k\f {2a+1}{2k+1}
\e 1+2t+p^2t(t+1)B_{p-2}(-a))\mod {p^3},\tag 1.7
\\&\q\sum_{k=1}^{p-1}\f 1k\b ak\b{-1-a}k
\e -\f 23p^2t(t+1)B_{p-3}(-a)
-2\f{B_{p^2(p-1)}(-a)-B_{p^2(p-1)}}{p^2(p-1)} \mod {p^3},\tag 1.8
\\&\q\sum_{k=0}^{p-1}\b ak(-2)^k\e (-1)^{\ap}-ptE_{p-2}(-a)\mod {p^2},
\tag 1.9
\\&\q\ \sum_{k=1}^{p-1}\f {(-1)^{k-1}}k\b ak\e
\f{B_{p^2(p-1)}(-a)-B_{p^2(p-1)}}{p^2(p-1)}
-\f{a-\ap}2B_{p-2}(-a)\mod{p^2}.\tag 1.10\endalign$$  We note that
Tauraso[T2] obtained a congruence for $\sum_{k=1}^{p-1}\f 1k\b
ak\b{-1-a}k\mod{p^2}$, and he proved that ([T1])
$$\sum_{k=1}^{p-1}\f{(-1)^{k-1}}k\b{-1/2}k\e H_{\f{p-1}2}\mod{p^3},$$
where $H_n=1+\f 12+\cdots +\f 1n$. Z.W. Sun proved that ([Su1])
$$\sum_{k=0}^{p-1}\f{\b{2k}k^2}{16^k}=\sum_{k=0}^{p-1}\b{-\f
12}k^2\e 1+(-1)^{\f{p-1}2}-p^2E_{p-3}\mod{p^3}.$$ Recently, using
the fact
$$\b{a-1}{\f{p-1}2}\b{-a-1}{\f{p-1}2}
\e\cases \f {pt}{\ap}\mod{p^2}&\t{if $\ap\le \f{p-1}2$,}
\\\f{p(t+1)}{\ap}\mod{p^2}&\t{if $\ap>\f{p-1}2$}
\endcases\tag 1.11$$
 and the
method in [S9-S11] Mao and Sun [MS] obtained congruences for $
\sum_{k=0}^{(p-1)/2}\b ak\b{-1-a}k$ and $\sum_{k=0}^{(p-1)/2}\b
ak\b{-1-a}k\f 1{2k+1}$ modulo $p^2$. In particular, they proved
$$\align&\sum_{k=0}^{(p-1)/2}\f{\b{2k}k\b{3k}k}{27^k}\e \Ls
p3\f{2^p+1}3\mod {p^2},\
\sum_{k=0}^{(p-1)/2}\f{\b{2k}k\b{4k}{2k}}{(2k+1)64^k}\e \Ls
{-1}p2^{p-1}\mod{p^2},
\\&\sum_{k=0}^{(p-1)/2}\f{\b{3k}k\b{6k}{3k}}{(2k+1)432^k}\e
\Ls p3\f{3^p+1}4\mod{p^2}.\endalign$$

\par Let $p>3$ be a prime and $a\in\Bbb Z_p$ with
$a\not\e 0\mod p$. In this paper we establish congruences for
$$\sum_{k=1}^{(p-1)/2}\f{\b ak\b{-1-a}k}k, \q\sum_{k=0}^{(p-1)/2}k\b
ak\b{-1-a}k \qtq{and}\sum_{k=0}^{(p-1)/2}\f{\b
ak\b{-1-a}k}{2k-1}\mod {p^2}$$ and some transformation formulas for
congruences modulo $p^2$. For instance, for $\ap<\f p2$ we have
$$\align&\sum_{k=1}^{(p-1)/2}\f{\b ak\b{-1-a}k}k\e -2\sum_{k=1}^{p-1}
\f{(-1)^{k-1}}k\b ak \mod{p^2} \\&\sum_{k=0}^{(p-1)/2}\b ak\b{-1-a}k
\e \sum_{k=0}^{p-1}\b ak(-2)^k\mod {p^2},
\\&\sum_{k=0}^{(p-1)/2}\f{\b ak\b{-1-a}k}{2k-1}\e
2(a-\ap)-(2a+1)\sum_{k=0}^{p-1}\b{2a}k(-2)^k\mod{p^2}.
\endalign$$
 As consequences and applications we obtain  some new congruences modulo $p^2$. Here
 are three typical examples:
$$\align &\sum_{k=0}^{(p-1)/2}\f{\b{2k}k\b{4k}{2k}}{(2k-1)64^k} \e
(-1)^{\f{p+1}2}\f{p+1}2\mod{p^2},\tag 1.12
\\&\sum_{k=1}^{(p-1)/2}\f{\b{2k}k\b{4k}{2k}}{k\cdot 64^k}
 \e
 6q_p(2)-p\big(3q_p(2)^2+2(-1)^{\f{p-1}2}E_{p-3}\big)\mod{p^2},\tag
 1.13
  \\&\sum_{k=0}^{(p-1)/2}\f{\b{2k}k\b{4k}{2k}}{64^k}
\e (-1)^{[\f p4]}\big(1+P_{p-\sls 2p}\big)\mod{p^2}.\tag 1.14
 \endalign$$
 where $q_p(a)=(a^{p-1}-1)/p$ and  $\{P_n\}$ is the Pell sequence given by
 $P_0=0,\ P_1=1$ and $P_{n+1}=2P_n+P_{n-1}\ (n\ge 1).$
 We note that (1.13) is equivalent to a conjecture made by the author's brother Z.W. Sun
 in [Su2].

 \section* {2. Congruences for $\sum_{k=1}^{\f{p-1}2}
 \f 1k\b ak\b{-1-a}k\mod {p^2}$}
\par\q For any positive integer $n$ and variable $a$ let
$$S_n(a)=\sum_{k=1}^n\f 1k\b ak\b{-1-a}k.$$
By [S9, (2.1)],
$$S_n(a)-S_n(a-1)=\f 2a\b{a-1}n \b{-a-1}n-\f 2a.\tag 2.1$$
Let $p>3$ be a prime, $a\in\Bbb Z_p$ and  $t=(a-\ap)/p$. Using (2.1)
we see that

$$\align &S_n(a)-S_n(a-\ap)
\\&=\sum_{k=0}^{\ap-1}(S_n(a-k)-S_n(a-k-1))
=\sum_{k=0}^{\ap-1}\f 2{a-k}\Big\{\b{a-k-1}n \b{k-a-1}n-1\Big\}
\\&=\sum_{k=0}^{\ap-1}\f 2{pt+\ap-k}
\Big\{\b{pt+\ap-k-1}n \b{-pt-(\ap-k)-1}n-1\Big\} .\endalign$$ Hence
$$S_n(a)-S_n(pt)=\sum_{r=1}^{\ap}\f 2{pt+r}
\Big\{\b{pt+r-1}n \b{-pt-r-1}n-1\Big\}.\tag 2.2$$
\par For any positive integer $n$ define $H_n=\sum_{k=1}^n\f 1k$.
For convenience we also define $H_0=0$.
 \pro{Lemma 2.1}
 Let $p>3$ be a prime, $r\in\{1,2,\ldots,p-1\}$ and $t\in\Bbb
Z_p$. Then
$$\align &\b{pt+r-1}{\f{p-1}2}\b{-pt-r-1}{\f{p-1}2}
\\&\e\cases
\f{pt}r+\f{p^2t}r\Big(2q_p(2)+H_{\f{p-1}2-r}\Big)-\f{p^2t^2}{r^2}\mod{p^3}
&\t{if $r<\f p2$,}
\\\f{p(t+1)}r+\f{p^2(t+1)}r\Big(2q_p(2)+H_{r-\f{p+1}2}\Big)-\f{p^2t(t+1)}{r^2}\mod{p^3}
&\t{if $r>\f p2$.}\endcases\endalign$$
\endpro
Proof. We first assume $r<\f p2$. It is clear that
$$\align &{\big(\f{p-1}2\big)!}^2
\b{pt+r-1}{\f{p-1}2}\b{-pt-r-1}{\f{p-1}2}
\\&={\big(\f{p-1}2\big)!}^2\b{pt+r-1}{\f{p-1}2}(-1)^{\f{p-1}2}\b{pt+r+\f{p-1}2}{\f{p-1}2}
\\&=(-1)^{\f{p-1}2}\f{pt}{pt+r}(pt-1)\cdots\Big(pt-\big(\f{p-1}2-r\big)\Big)(pt+1)(pt+2)\cdots
\Big(pt+\f{p-1}2+r\Big)
\\&=(-1)^{\f{p-1}2}\f{pt}{pt+r}((pt)^2-1^2)((pt)^2-2^2)\cdots
\Big((pt)^2-\big(\f{p-1}2-r\big)^2\Big)
\\&\q\times\prod_{s=0}^{r-1}\Big(pt+\f{p-1}2-s\Big)
\Big(pt+\f{p-1}2+s+1\Big)
\\&\e (-1)^r\f{pt}{pt+r}{\Big(\f{p-1}2-r\Big)!}^2
\ \prod_{s=0}^{r-1}\Big(\f{p-1}2-s\Big) \Big(\f{p-1}2+s+1\Big)
\\&=(-1)^r\f{pt}{pt+r}\Big(\f{p-1}2-r\Big)!\Big(\f{p-1}2+r\Big)!
=(-1)^r\f{pt}{pt+r}\cdot\f{(p-1)!}{\b{p-1}{\f{p-1}2-r}}
\\&\e
(-1)^r\f{pt}{pt+r}\cdot\f{(p-1)!}{(-1)^{\f{p-1}2-r}(1-pH_{\f{p-1}2-r})}
\e (-1)^{\f{p-1}2}(p-1)!\cdot pt\f{(r-pt)(1+pH_{\f{p-1}2-r})}{r^2}
\\&\e (-1)^{\f{p-1}2}(p-1)!\cdot pt\Big(\f 1r-\f{pt}{r^2}+\f{pH_{\f{p-1}2-r}}r\Big)
\mod{p^3}.
\endalign$$
Since
$$\f{(p-1)!}{{(\f{p-1}2)!}^2}=\f{(p-1)\cdots(p-\f{p-1}2)}{(\f{p-1}2)!}
\e (-1)^{\f{p-1}2}\big(1-pH_{\f{p-1}2}\big)\e
(-1)^{\f{p-1}2}(1+2pq_p(2))\mod {p^2},$$ from the above we deduce
that
$$\align\b{pt+r-1}{\f{p-1}2}\b{-pt-r-1}{\f{p-1}2}&\e (1+2pq_p(2))pt
\Big(\f 1r-\f{pt}{r^2}+\f{pH_{\f{p-1}2-r}}r\Big)
\\&\e pt\Big(\f 1r+\f{p(2q_p(2)+H_{\f{p-1}2-r})}r-\f{pt}{r^2}\Big)\mod
{p^3}.\endalign$$ This yields the result in the case $r<\f p2$.
\par Now assume $r>\f p2$. Set $t'=-t-1$ and $r'=p-r$. We see that
$$\align \b{pt+r-1}{\f{p-1}2}\b{-pt-r-1}{\f{p-1}2}
&=\b{p(t+1)-(p-r)-1}{\f{p-1}2}\b{-p(t+1)+p-r-1}{\f{p-1}2}
\\&=\b{pt'+r'-1}{\f{p-1}2}\b{-pt'-r'-1}{\f{p-1}2}.\endalign$$
Since $r'<\f p2$, from the above we deduce that
$$\align &\b{pt+r-1}{\f{p-1}2}\b{-pt-r-1}{\f{p-1}2}
\\&=\b{pt'+r'-1}{\f{p-1}2}\b{-pt'-r'-1}{\f{p-1}2}
\e \f{pt'}{r'}+\f{p^2t'}{r'}\Big(2q_p(2)+H_{\f{p-1}2-r'}\Big)
-\f{p^2{t'}^2}{{r'}^2}
\\&=\f{p(t+1)}{r-p}-\f{p^2(t+1)}{p-r}\Big(2q_p(2)+H_{r-\f{p+1}2}\Big)-\f{p^2(t+1)^2}
{(p-r)^2}
\\&\e \f{p(t+1)(r+p)}{r^2}+\f{p^2(t+1)}r\Big(2q_p(2)+H_{r-\f{p+1}2}\Big)
-\f{p^2(t+1)^2}{r^2} \mod{p^3}.\endalign$$ This yields the result in
the case $r>\f p2$. Hence the lemma is proved.
\par We remark that Lemma 2.1 improves Mao and Sun's (1.11).

 \pro{Theorem 2.1} Let $p>3$ be a prime, $a\in\Bbb Z_p$ and
$a\not\e 0\mod p$. Then
$$ \sum_{k=1}^{(p-1)/2}\f{\b ak\b{-1-a}k}k+2\sum_{k=1}^{p-1}
\f{(-1)^{k-1}}k\b ak\e\cases 0\mod {p^2}&\t{if $\ap\le \f{p-1}2$,}
\\pB_{p-2}(-a)\mod{p^2}&\t{if $\ap>\f{p-1}2$}
\endcases$$
and
$$\align &\sum_{k=1}^{p-1}
\f{(-1)^{k-1}}k\b ak \\&\e
-\f{B_{2p-2}(-a)-B_{2p-2}}{2p-2}+2\f{B_{p-1}(-a)-B_{p-1}}{p-1}-\f
{a-\ap}2B_{p-2}(-a)\mod {p^2}.\endalign$$
\endpro
Proof. Set $t=(a-\ap)/p$. By [S3, Theorem 5.2],
$\sum_{k=1}^{(p-1)/2}\f 1{k^2}\e 0\mod p$. Thus,
$$\align S_{\f{p-1}2}(pt)&=\sum_{k=1}^{(p-1)/2}\f
1k\b{pt}k\b{-pt-1}k =\sum_{k=1}^{(p-1)/2}\f{pt}{k^2}\b{pt-1}{k-1}
\b{-pt-1}k
\\&\e \sum_{k=1}^{(p-1)/2}\f{pt}{k^2}
\b{-1}{k-1}\b{-1}k =-pt\sum_{k=1}^{(p-1)/2}\f 1{k^2}\e
0\mod{p^2}.\endalign$$ For $1\le \ap\le \f{p-1}2$, from (2.2) and
Lemma 2.1 we see that
$$\align S_{\f{p-1}2}(a)&=S_{\f{p-1}2}(pt)
+\sum_{r=1}^{\ap}\f{2(r-pt)}{r^2-p^2t^2}
\Big\{\b{pt+r-1}{\f{p-1}2}\b{-pt-r-1}{\f{p-1}2}-1\Big\}
\\&\e 2\sum_{r=1}^{\ap}\f{r-pt}{r^2}\Big(\f{pt}r-1\Big)
\e 2\sum_{r=1}^{\ap}\Big( \f{2pt}{r^2}-\f 1r\Big)
\mod{p^2}.\endalign$$ For $\ap>\f{p-1}2$, from (2.2) and Lemma 2.1
we deduce that
$$\align S_{\f{p-1}2}(a)&=S_{\f{p-1}2}(pt)
+\sum_{r=1}^{\ap}\f{2(r-pt)}{r^2-p^2t^2}
\Big\{\b{pt+r-1}{\f{p-1}2}\b{-pt-r-1}{\f{p-1}2}-1\Big\}
\\&\e 2\sum_{r=1}^{(p-1)/2}
\f{r-pt}{r^2}\Big(\f{pt}r-1\Big)+2\sum_{r=(p+1)/2}^{\ap}
\f{r-pt}{r^2}\Big(\f{p(t+1)}r-1\Big)
\\&=2\sum_{r=1}^{\ap}\f{r-pt}{r^2}\Big(\f{p(t+1)}r-1\Big)
-2\sum_{r=1}^{(p-1)/2}\f{r-pt}{r^2}\cdot \f pr
\\&\e 2\sum_{r=1}^{\ap}\Big( \f{2pt}{r^2}-\f 1r\Big)+2p\sum_{r=1}^{\ap}
\f 1{r^2}-2p\sum_{r=1}^{(p-1)/2}\f 1{r^2}\mod {p^2}.
\endalign$$
By [S3, Theorem 5.2],  $\sum_{r=1}^{(p-1)/2}\f 1{r^2}\e 0\mod p$. By
[S9, (3.4)] and Fermat's little theorem, $\sum_{r=1}^{\ap}\f
1{r^2}\e \f 12B_{p-2}(-a)\mod p$. From the above we deduce that
$$ S_{\f{p-1}2}(a)-2\sum_{r=1}^{\ap}\Big( \f{2pt}{r^2}-\f 1r\Big)
\e\cases 0\mod {p^2}&\t{if $\ap\le \f{p-1}2$,}
\\pB_{p-2}(-a)\mod{p^2}&\t{if $\ap>\f{p-1}2$.}
\endcases\tag 2.3$$
  Putting $n=p-1$ and $b=-1$ in [S9, Lemma 3.1] and
then applying [S3, Theorem 5.1] we get
$$\align \sum_{k=1}^{p-1}\f{(-1)^k}k\b{a}k
&\e -pt\sum_{k=1}^{p-1}\f 1{k^2}-\sum_{r=1}^{\ap} \f
1r+2pt\sum_{r=1}^{\ap}\f 1{r^2} \e \sum_{r=1}^{\ap}\Big(
\f{2pt}{r^2}-\f 1r\Big) \mod{p^2}.\endalign$$ This together with
(2.3) yields the the first part.
\par Putting $b=p-1$ and $x=-a$ in [S3, Theorem 3.1] and then
applying [S3, Theorem 2.2 (with $k=p^2-1$ and $n=2$)] we deduce that
$$\align \f{B_{p^2(p-1)}(-a)-B_{p^2(p-1)}}{p^2(p-1)}
&=\f{B_{(p^2-1)(p-1)+p-1}(-a)-B_{(p^2-1)(p-1)+p-1}}{(p^2-1)(p-1)+p-1}
\\&\e (p^2-1)\f{B_{p-1+p-1}(-a)-B_{p-1+p-1}}{p-1+p-1}-(p^2-2)
\f{B_{p-1}(-a)-B_{p-1}}{p-1}
\\&\e -\f{B_{2p-2}(-a)-B_{2p-2}}{2p-2}+2\f{B_{p-1}(-a)-B_{p-1}}{p-1}
\mod{p^2}.\endalign$$ This together with (1.10) ([S9, Theorem 3.1])
yields the remaining part. The proof is now complete.

 \pro{Theorem 2.2} Let $p>3$ be a prime. Then
 $$\align &\sum_{k=1}^{(p-1)/2}\f{\b{2k}k\b{4k}{2k}}{k\cdot 64^k}
 \e
 6q_p(2)-p\big(3q_p(2)^2+2(-1)^{\f{p-1}2}E_{p-3}\big)\mod{p^2},\tag i
 \\&\sum_{k=1}^{(p-1)/2}\f{\b{2k}k\b{3k}{k}}{k\cdot 27^k}
\e 3q_p(3)-p\Big(\f 32q_p(3)^2 +2\Ls p3U_{p-3}\Big)\mod{p^2},\tag ii
 \\&\sum_{k=1}^{(p-1)/2}\f{\b{3k}k\b{6k}{3k}}{k\cdot 432^k}
 \e  4q_p(2)+3q_p(3)-p\Big(2q_p(2)^2+\f 32q_p(3)^2+5\Ls
 p3U_{p-3}\Big)
 \mod{p^2}.
 \tag iii
 \endalign$$
 \endpro
Proof. From [S4, Lemma 2.5] we know that
$E_{2n}=-4^{2n+1}\f{B_{2n+1}\sls 14}{2n+1}$. Thus,
$E_{p-3}=-4^{p-2}\f{B_{p-2}\sls 14}{p-2}$ $\e \f 18B_{p-2}\sls
14\mod p$. Now taking $a=-\f 14$ in Theorem 2.1 we see that
$$\align &\sum_{k=1}^{(p-1)/2}\f{\b{-1/4}k\b{-3/4}k}k+2
\sum_{k=1}^{p-1}\f{(-1)^{k-1}}k\b{-1/4}k\\&
\e\cases 0\mod p &\t{if
$p\e 1\mod 4$,}
\\ pB_{p-2}\ls 14\e 8pE_{p-3}\mod p&\t{if $p\e 3\mod 4$.}
\endcases\endalign$$
By [S9, Theorem 3.2],
$$-2\sum_{k=1}^{p-1}\f{(-1)^{k-1}}k\b{-1/4}k \e
6q_p(2)+p(-3q_p(2)^2-2(2-(-1)^{\f{p-1}2}) E_{p-3})\mod{p^2}.$$ Hence
(i) is true by (1.1).
 \par Taking $a=-\f 13$ in Theorem 2.1 and  applying the fact $B_{p-2}(\f
13)\e 6U_{p-3}\mod p$ ([S6]) we see that
$$\align &\sum_{k=1}^{(p-1)/2}\f{\b{-1/3}k\b{-2/3}k}k+2
\sum_{k=1}^{p-1}\f{(-1)^{k-1}}k\b{-1/3}k\\& \e\cases 0\mod p &\t{if
$p\e 1\mod 3$,}
\\ pB_{p-2}\ls 13\e 6pU_{p-3}\mod p&\t{if $p\e 2\mod 3$.}
\endcases\endalign$$
By [S9, Theorem 3.3],
$$-2\sum_{k=1}^{p-1}\f{(-1)^{k-1}}k\b{-1/3}k \e
3q_p(3)-p\Big(\f 32q_p(3)^2+\big(3-\sls p3\big)
U_{p-3}\Big)\mod{p^2}.$$ Thus (ii) holds by (1.1).

\par Taking $a=-\f 16$ in Theorem 2.1 and  applying the fact $B_{p-2}(\f
16)\e 30U_{p-3}\mod p$ ([S6]) we see that
$$\align &\sum_{k=1}^{(p-1)/2}\f{\b{-1/6}k\b{-5/6}k}k+2
\sum_{k=1}^{p-1}\f{(-1)^{k-1}}k\b{-1/6}k\\& \e\cases 0\mod p &\t{if
$p\e 1\mod 3$,}
\\ pB_{p-2}\ls 16\e 30pU_{p-3}\mod p&\t{if $p\e 2\mod 3$.}
\endcases\endalign$$
By [S9, Theorem 3.4],
$$\align &-2\sum_{k=1}^{p-1}\f{(-1)^{k-1}}k\b{-1/6}k
\\& \e 4q_p(2)+3q_p(3)-p\Big(2q_p(2)^2+\f 32q_p(3)^2+5
\big(3-2\sls p3\big) U_{p-3}\Big)\mod{p^2}.\endalign$$ Thus (iii)
holds by (1.1).

\section*{3. Congruences for $\sum_{k=0}^{\f{p-1}2}
\b ak\b{-1-a}k$ and $\sum_{k=0}^{\f{p-1}2}k \b
ak\b{-1-a}k\mod{p^2}$}
\par Let $p$ be a prime greater than $3$, $a\in\Bbb Z_p$ and $t=(a-\ap)/p$. From [S11, p.3299] we know that for any positive integer $n$,
$$\sum_{k=0}^n\b ak\b{-1-a}k-(-1)^{\ap}\sum_{k=0}^n\b{pt}k\b{-1-pt}k
=2\sum_{k=0}^{\ap-1}(-1)^k\b{a-k-1}n\b{k-a-1}n.$$ For $n<p$ we see
that
$$\align \sum_{k=0}^n\b{pt}k\b{-1-pt}k
&=1+\sum_{k=1}^n\f{pt}k\b{pt-1}{k-1} \b{-1-pt}k\\&\e
1+pt\sum_{k=1}^n\f 1k\b{-1}{k-1}\b{-1}k=1-ptH_n\mod{p^2}.\endalign$$
Taking $n=\f{p-1}2$ in the above and then applying (1.11) or Lemma
2.1 we deduce the following lemma due to Mao and Sun.
 \pro{Lemma 3.1 ([MS, (3.5)]} Let $p>3$ be a prime,
$a\in\Bbb Z_p$, $a\not\e 0\mod p$ and $t=(a-\ap)/p$. Then
$$\sum_{k=0}^{(p-1)/2}\b ak\b{-1-a}k\e (-1)^{\ap}(1-ptH_{\f{p-1}2})+2p
\sum_{k=0}^{\ap-1}(-1)^{k}\f{t+\delta_k}{\ap-k}\mod{p^2},$$ where
$\delta_k=1$ or $0$ according as $\ap-k>\f{p-1}2$ or not.
\endpro
\pro{Theorem 3.1} Let $p>3$ be a prime, $a\in\Bbb Z_p$, $a\not\e
0\mod p$ and $t=(a-\ap)/p$. Then
$$\align&\sum_{k=0}^{(p-1)/2}\b ak\b{-1-a}k
\\&\e\cases \sum_{k=0}^{p-1}\b ak(-2)^k\mod {p^2}&\t{if $\ap\le
\f{p-1}2$,}
\\\big(1+\f 1t\big)\sum_{k=0}^{p-1}\b ak(-2)^k-(-1)^{\ap}\f 1t
\mod {p^2}&\t{if $\ap>\f{p-1}2$ and $t\not\e 0\mod p$,}
\\(-1)^{\ap}-pE_{p-2}(-a)\mod {p^2}&\t{if $\ap>\f{p-1}2$ and $t\e 0\mod p$.}
\endcases\endalign$$
\endpro
Proof. Note that
$$\aligned \sum_{k=0}^{\ap-1}(-1)^{\ap-k}\f{t+\delta_k}{\ap-k}
&=\sum\Sb 0\le k\le \ap-1\\\ap-k<\f p2\endSb (-1)^{\ap-k}\f t{\ap-k}
+\sum\Sb 0\le k\le \ap-1\\\ap-k>\f
p2\endSb(-1)^{\ap-k}\f{t+1}{\ap-k}
\\&=\sum\Sb 1\le r\le \ap\\r<\f p2\endSb(-1)^r\f tr+
\sum\Sb 1\le r\le \ap\\r>\f p2\endSb(-1)^r\f {t+1}r
\\&=\cases \sum_{r=1}^{\ap}(-1)^r\f tr&\t{if $\ap\le\f{p-1}2$,}
\\\sum_{r=1}^{\ap}(-1)^r\f {t+1}r-\sum_{r=1}^{(p-1)/2}
\f{(-1)^r}r&\t{if $\ap>\f{p-1}2$.}
\endcases\endaligned$$
From [L] we know that $H_{\f{p-1}2}\e -2q_p(2)\mod p$ and $H_{[\f
p4]}\e -3q_p(2)\mod p$. Thus,
$$\sum_{r=1}^{(p-1)/2}\f{(-1)^r}r
=\sum_{r=1}^{(p-1)/2}\f{(-1)^r+1}r-H_{\f{p-1}2} = H_{[\f
p4]}-H_{\f{p-1}2} \e -q_p(2)\mod p.$$ Hence applying Lemma 3.1 we
obtain
$$\align\sum_{k=0}^{(p-1)/2}\b ak\b{-1-a}k
&\e (-1)^{\ap}(1+2ptq_p(2))+2pt(-1)^{\ap}\sum_{r=1}^{\ap}\f{(-1)^r}r
\\&\q+\cases
0\mod{p^2}&\t{if $\ap\le \f{p-1}2$,}
\\2p(-1)^{\ap}\sum_{r=1}^{\ap}\f{(-1)^r}r+2p(-1)^{\ap}q_p(2)\mod{p^2}
 &\t{if $\ap>\f{p-1}2$.}
\endcases\endalign$$
By [S11, p.3306],
$$\sum_{r=1}^{\ap}\f{(-1)^r}r+q_p(2)\e -\f 12 (-1)^{\ap}E_{p-2}(-a)
\mod p.$$ Thus,
$$\aligned \sum_{k=0}^{\f{p-1}2}\b ak\b{-1-a}k
\e\cases
(-1)^{\ap}(1+2pt\big(\sum_{r=1}^{\ap}\f{(-1)^r}r+q_p(2)\big)
\\\qq\e
(-1)^{\ap}-ptE_{p-2}(-a)\mod {p^2}\qq\q\  \t{if $\ap\le \f{p-1}2$,}
\\(-1)^{\ap}(1+2p(t+1)\big(\sum_{r=1}^{\ap}\f{(-1)^r}r+q_p(2)\big)
\\\qq\e
(-1)^{\ap}-p(t+1)E_{p-2}(-a)\mod {p^2}\q \t{if $\ap>
\f{p-1}2$.}\endcases\endaligned\tag 3.1$$ By [S11, Theorem 3.1],
$$ (-1)^{\ap}-ptE_{p-2}(-a)
\e \sum_{k=0}^{p-1}\b ak(-2)^k\mod {p^2}.$$ Thus,  the result is
true for $\ap\le\f{p-1}2$. For $\ap>\f{p-1}2$ and $t\not\e 0\mod p$
we have
$$\align \sum_{k=0}^{(p-1)/2}\b ak\b{-1-a}k
&\e (-1)^{\ap}+\f{t+1}t(-ptE_{p-2}(-a))
\\&\e (-1)^{\ap}+\f{t+1}t\Big(\sum_{k=0}^{p-1}\b ak(-2)^k-(-1)^{\ap}\Big)
\mod{p^2}.\endalign$$ This yields the result in this case. The proof
is now complete.

\pro{Theorem 3.2} Let $p>3$ be a prime. Then
$$\sum_{k=0}^{(p-1)/2}\f{\b{2k}k\b{4k}{2k}}{64^k}
\e (-1)^{[\f p4]}\big(1+P_{p-\sls 2p}\big)\mod{p^2}.$$
\endpro
Proof. Clearly
$$\align \sum_{r=1}^{\langle-1/4\rangle_p}\f{(-1)^r}r
&=\sum_{r=1}^{\langle-1/4\rangle_p}\f{(-1)^r+1}r-\sum_{r=1}^{\langle-1/4\rangle_p}\f
1r
 =\cases H_{[\f p8]}-H_{[\f p4]}&\t{if $p\e 1\mod 4$,}
 \\H_{[\f {3p}8]}-H_{[\f {3p}4]}&\t{if $p\e 3\mod 4$.}
 \endcases\endalign$$
 From [S1, Theorems 3.3 and 3.4] we know that
 $$\align &H_{[\f p8]}\e -4q_p(2)-2\f{P_{p-\sls 2p}}p\mod p,\
\\&H_{[\f {3p}8]}=H_{[\f p8]}+\sum_{\f p8<k<\f{3p}8}\f 1k\e H_{[\f p8]}
+4\f{P_{p-\sls 2p}}p\e  -4q_p(2)+2\f{P_{p-\sls 2p}}p\mod p.
\endalign$$
We also have $H_{[\f{3p}4]}=H_{p-1}-\sum_{k=1}^{[p/4]}\f 1{p-k} \e
H_{[\f p4]}\e -3q_p(2)\mod p$. Thus,
$$\align\sum_{r=1}^{\langle-1/4\rangle_p}\f{(-1)^r}r
 =\cases H_{[\f p8]}-H_{[\f p4]}\e -q_p(2)-2 \f{P_{p-\sls 2p}}p\mod
p&\t{if $p\e 1\mod 4$,}
 \\H_{[\f {3p}8]}-H_{[\f {3p}4]}\e -q_p(2)
 +2\f{P_{p-\sls 2p}}p\mod p&\t{if $p\e 3\mod 4$.}
 \endcases\endalign$$
 Taking $a=-1/4$ in (3.1) we see that $t=(
(-1)^{\f{p-1}2}-2)/4$ and so
$$\align&\sum_{k=0}^{(p-1)/2}\f{\b{2k}k\b{4k}{2k}}{64^k}
\\&=\sum_{k=0}^{(p-1)/2}\b{-1/4}k\b{-3/4}k
 \\&\e \cases(-1)^{\ap}(1+2pt(-2)\f{P_{p-\sls
 2p}}p)=(-1)^{\f{p-1}4}(1+P_{p-\sls 2p})\mod p&\t{if $p\e 1\mod 4$,}
 \\(-1)^{\ap}(1+2p(t+1)\cdot 2\f{P_{p-\sls
 2p}}p)=(-1)^{\f{p-3}4}(1+P_{p-\sls 2p})\mod p&\t{if $p\e 3\mod 4$.}
 \endcases\endalign$$
 This proves the theorem.

 \par It is well known that ([MOS]) for any positive integer $n$,
 $$\aligned &B_{2n}\Ls 12=(2^{1-2n}-1)B_{2n},\ B_{2n}\Ls 13=\f{3-3^{2n}}{2\cdot 3^{2n}}B_{2n},\q
 \\&B_{2n}\Ls 14=\f{2-2^{2n}}{4^{2n}}B_{2n},\
 B_{2n}\Ls 16=\f{(2-2^{2n})(3-3^{2n})}{2\cdot 6^{2n}}B_{2n}.
 \endaligned\tag 3.2$$
 Let $p>3$ be a prime. From (3.2) and the well-known fact
$pB_{p-1}\e p-1\mod p$ (see [IR]) one can easily deduce the
 following known congruences (see [L],[GS]):
 $$\align&B_{p-1}\Ls 12-B_{p-1}\e 2q_p(2)\mod p,
 \q B_{p-1}\Ls 13-B_{p-1}\e\f 32q_p(3)\mod p,\tag 3.3
 \\&B_{p-1}\Ls 14-B_{p-1}\e 3q_p(2)\mod p,\tag 3.4
 \\&B_{p-1}\Ls 16-B_{p-1}\e 2q_p(2)+\f 32q_p(3)\mod p.\tag 3.5
 \endalign$$
 In [GS] Graville and Sun showed that
 $$B_{p-1}\Ls 1{12}-B_{p-1}\e 3\f{S_{p-\sls 3p}}p+3q_p(2)+\f
32q_p(3)\mod p,\tag 3.6$$ where $\{S_n\}$ is given by $S_0=0,\
S_1=1$ and $S_{n+1}=4S_n-S_{n-1}$ $(n\ge 1)$.

 \pro{Theorem 3.3} Let $p>3$ be a prime. Then
$$\sum_{k=0}^{(p-1)/2}\f{\b{3k}k\b{6k}{3k}}{432^k}
\e (-1)^{\f {p-1}2}\Big(1+\Ls 3pS_{p-\sls 3p}\Big)\mod{p^2}.$$
\endpro
Proof. It is well known that ([MOS])
$$E_n(x)=\f 2{n+1}\Big(B_{n+1}(x)-2^{n+1}B_{n+1}\Ls x2\Big)
=\f{2^{n+1}}{n+1}\Big(B_{n+1}\Ls{x+1}2-B_n\Ls x2\Big).\tag 3.7$$ For
$a\in\Bbb Z_p$ we know that $B_{p-1}(a)-B_{p-1}\in\Bbb Z_p$ and
$pB_{p-1}\e p-1\mod p$ (see [IR,S2,S3]). Thus,
$$\aligned E_{p-2}(-a)&=\f 2{p-1}\Big(B_{p-1}(-a)-2^{p-1}B_{p-1}
\Big(-\f a2\Big)\Big)
\\&=\f 2{p-1}\Big(B_{p-1}(-a)-B_{p-1}-2^{p-1}
\Big(B_{p-1} \Big(-\f
a2\Big)-B_{p-1}\Big)-\big(2^{p-1}-1\big)B_{p-1}\Big)
\\&\e -2\Big(B_{p-1}(-a)-B_{p-1}-\Big(B_{p-1} \Big(-\f
a2\Big)-B_{p-1}\Big)+q_p(2)\Big)\mod p.\endaligned\tag 3.8$$ Taking
$a=-\f 16$ in (3.8) and then applying (3.5) and (3.6) we get
$$E_{p-2}\Ls 16\e -2\Big(B_{p-1}\Ls 16-B_{p-1}-\Big(B_{p-1}\Ls 1{12}
-B_{p-1} \Big)+q_p(2)\Big)\e 6\f{S_{p-\sls 3p}}p\mod p.\tag 3.9$$
Now applying (3.1) we deduce that
$$\align
\sum_{k=0}^{(p-1)/2}\b{-\f 16}k\b{-\f 56}k  \e\cases
(-1)^{\f{p-1}6}+\f p6E_{p-2}\sls 16\e (-1)^{\f{p-1}2}+S_{p-\sls
3p}\mod {p^2}&\t{if $6\mid p-1$,}
\\(-1)^{\f{5p-1}6}-\f p6E_{p-2}\sls 16\e (-1)^{\f{p-1}2}-S_{p-\sls
3p}\mod {p^2}&\t{if $6\mid p-5$.}
\endcases\endalign$$
This together with (1.1) yields the result.
\pro{Theorem 3.4} Let
$p>3$ be a prime, $a\in\Bbb Z_p$, $a\not\e 0\mod p$ and
$t=(a-\ap)/p$. Then
$$\align &\sum_{k=0}^{(p-1)/2}k\b ak\b{-1-a}k
\\&\e\cases (-1)^{\ap}a(a+1)-\f 12pt(2a(a+1)E_{p-2}(-a)+2a+1)\mod
{p^2}&\t{if $\ap\le \f{p-1}2$,}
\\(-1)^{\ap}a(a+1)-\f 12p(t+1)(2a(a+1)E_{p-2}(-a)+2a+1)\mod
{p^2}&\t{if $\ap>\f{p-1}2$.}
\endcases\endalign$$\endpro
Proof. By [S11, Lemma 2.4],
$$\align &\sum_{k=0}^n(k-a(a+1))\b ak\b{-1-a}k
\\&=-a(a+1)\b{a-1}n\b{-2-a}n =-a(a+n+1)\b{a-1}n\b{-1-a}n.\endalign$$
Taking $n=\f{p-1}2$ in the above identity and then applying Lemma
2.1 and (3.1) we deduce that
$$\align &\sum_{k=0}^{(p-1)/2}k\b ak\b{-1-a}k
- a(a+1)\sum_{k=0}^{(p-1)/2}\b ak\b{-1-a}k
\\&=-a\Big(a+\f{p-1}2+1\Big)\b{a-1}{\f{p-1}2}\b{-1-a}{\f{p-1}2}
\\&\e \cases -a(a+\f 12)\f{pt}{\ap}\e -\f{2a+1}2pt\mod{p^2}&\t{if $\ap\le
\f{p-1}2$,}
\\-a(a+\f 12)\f{p(t+1)}{\ap}\e -\f{2a+1}2p(t+1)\mod{p^2}&\t{if $\ap>\f{p-1}2$.}
\endcases\endalign$$
This together with (3.1) yields the result.
 \pro{Theorem 3.5} Let $p$ be a prime greater than $3$. Then
 $$\sum_{k=0}^{(p-1)/2}\f{k\b{2k}k\b{3k}k}{27^k}\e
 -\Ls p3\f{4-3p+2^{p+2}}{54}  \mod {p^2}.$$
 \endpro
 Proof. From  (3.3), (3.5) and (3.8)  we derive that
 $$\aligned E_{p-2}\Ls 13&\e-2\Big(B_{p-1}\Ls 13-B_{p-1}-\Big(B_{p-1}\Ls
 16-B_{p-1}\Big)+q_p(2)\Big)
 \\&\e -2\Big(\f 32q_p(3)-2q_p(2)-\f 32q_p(3)+q_p(2)\Big)=2q_p(2)\mod
 p. \endaligned\tag 3.10$$
 Now taking $a=-\f 13$ in Theorem 3.4 and then applying (1.1) and
 (3.10) we deduce the result.

\section*{4. Congruences for $\sum_{k=0}^{\f{p-1}2}\f 1{2k+1}\b
ak\b{-1-a}k$ and $\sum_{k=0}^{\f{p-1}2}\f 1{2k-1}\b ak\b{-1-a}k$
$\mod {p^2}$} \pro{Lemma 4.1 ([MS])} Let $p>3$ be a prime, $a\in\Bbb
Z_p$, $a\not\e 0\mod p$ and $t=(a-\ap)/p$. Then
$$\sum_{k=0}^{(p-1)/2}\b ak\b{-1-a}k\f{2a+1}{2k+1}
\e 1+2t+4ptq_p(2)+2pt\sum_{r=1}^{\ap}\f 1r+2p\sum\Sb 1\le r\le \ap
\\r>\f p2\endSb\f 1r\mod{p^2}.$$
\endpro
\pro{Theorem 4.1 } Let $p>3$ be a prime, $a\in\Bbb Z_p$, $a\not\e
0\mod p$ and $t=(a-\ap)/p$. Then
$$\align &\sum_{k=0}^{(p-1)/2}\b ak\b{-1-a}k\f{2a+1}{2k+1}
\\&\e \cases 1+2t+4ptq_p(2)-2pt(B_{p-1}(-a)-B_{p-1})\mod {p^2}&\t{if
$\ap\le \f{p-1}2$,}
\\1+2t+4p(t+1)q_p(2)-2p(t+1)(B_{p-1}(-a)-B_{p-1})\mod {p^2}&\t{if
$\ap>\f{p-1}2$.}\endcases\endalign$$
\endpro
Proof. By [S3, Lemma 3.2],
$$\aligned\sum_{r=1}^{\ap}\f 1r&\e \sum_{r=1}^{\ap}r^{p-2}
\e (-1)^{p-1}\f{B_{p-1}(-a)-B_{p-1}}{p-1}+(-a+\ap)B_{p-2} \\& \e
-(B_{p-1}(-a)-B_{p-1})\mod p.\endaligned\tag 4.1$$ Thus, for $\ap\le
\f{p-1}2$ the result follows from Lemma 4.1.
\par Now assume $\ap>\f{p-1}2$. As $H_{\f{p-1}2}\e -2q_p(2)\mod p$,
we see that
$$\sum\Sb 1\le r\le \ap
\\r>\f p2\endSb\f 1r=\sum_{r=1}^{\ap}\f 1r-\sum_{r=1}^{(p-1)/2}\f
1r\e \sum_{r=1}^{\ap}\f 1r+2q_p(2)\mod p.$$ Hence the result follows
from the above and Lemma 4.1.
\pro{Corollary 4.1} Let $p>3$ be a
prime. Then
$$\sum_{k=0}^{(p-1)/2}\f{\b{2k}k\b{3k}k}{(2k+1)27^k}
\e \Ls p3(2-2^{p+1}+3^p)\mod{p^2}.$$
\endpro
Proof. As $t=-\f 13$ or $-\f 23$ according as $p\e 1\mod 3$ or $p\e
2\mod 3$, taking $a=-\f 13$ in Theorem 4.1 and then applying (3.3)
yields
$$\sum_{k=0}^{(p-1)/2}\f 1{2k+1}\b{-\f 13}k\b{-\f 23}k
\e \Ls p3(1-4pq_p(2)+3pq_p(3))=\Ls p3(2-2^{p+1}+3^p)\mod p.$$ This
together with (1.1) gives the result.

\pro{Theorem 4.2 } Let $p>3$ be a prime, $a\in\Bbb Z_p$ and $a\not\e
0\mod p$.  Then
$$\align &\sum_{k=0}^{(p-1)/2}\b ak\b{-1-a}k\f 1{2k-1}
\\&\e\cases -(2a+1)+2(a-\ap)(1+(2a+1)E_{p-2}(-2a))\mod{p^2}
&\t{if $\ap<\f p2$,}
\\2a+1+2(p+a-\ap)(1+(2a+1)E_{p-2}(-2a))\mod{p^2}
&\t{if $\ap>\f p2$.}\endcases
\endalign$$
\endpro
Proof.  By [S10, Lemma 3.1],
$$\aligned &\sum_{k=0}^n\b ak\b{-1-a}k\f{(2a(a+1)+1)k-a(a+1)}{4k^2-1}
\\&=\f{a(a+1)}{2n+1}\b{a-1}n\b{-2-a}n=\f{a(a+n+1)}{2n+1}\b{a-1}n\b{-a-1}n.
\endaligned\tag 4.2$$
Set $t=(a-\ap)/p$. By Lemma 2.1,
$$\align &\f 1p\b{a-1}{\f{p-1}2}\b{-a-1}{\f{p-1}2}
\\&\e \cases
\f{t}{\ap}+\f{pt}a(2q_p(2)+H_{\f{p-1}2-\ap})-\f{pt^2}{a^2} \mod
{p^2}&\t{if $\ap<\f p2$,}
\\\f{t+1}{\ap}+\f{p(t+1)}a(2q_p(2)+H_{\ap-\f{p+1}2})-\f{pt(t+1)}{a^2} \mod
{p^2}&\t{if $\ap>\f p2$.}\endcases\endalign$$ For $\ap<\f p2$, using
(4.1) and the fact $B_n(1-x)=(-1)^nB_n(x)$ ([MOS]) we see that
$$H_{\f{p-1}2-\ap}\e -\Big(B_{p-1}\big(\f 12+a\big)-B_{p-1}\Big)
=-\Big(B_{p-1}\big(\f 12-a\big)-B_{p-1}\Big)\mod p.$$ For $\ap>\f
p2$, using (4.1) we see that
$$H_{\ap-\f{p+1}2}\e -\Big(B_{p-1}\big(\f 12-a\big)-B_{p-1}\Big)
\mod p.$$
 Thus, taking $n=\f{p-1}2$ in (4.2) and applying the above we obtain
$$\align &\sum_{k=0}^{(p-1)/2}\b ak\b{-1-a}k\f{(2a(a+1)+1)k-a(a+1)}{4k^2-1}
\\&\e \cases (a+\f{p+1}2)\Big(
\f{t(pt+\ap)}{\ap}+pt\big(2q_p(2)+B_{p-1}-B_{p-1}\big(\f
12-a\big)\big)-\f{pt^2}{a} \Big)\mod
{p^2}\\\qq\qq\qq\qq\qq\qq\qq\qq\t{if $\ap<\f p2$,}
\\(a+\f{p+1}2)\Big(\f{(t+1)(pt+\ap)}{\ap}+p(t+1)
\big(2q_p(2)+B_{p-1}-B_{p-1}\big(\f 12-a\big))\\\
\qq\q\qq-\f{pt(t+1)} {a}\Big) \mod {p^2}\q\t{if $\ap>\f
p2$}\endcases
\\&\e t'\Big(\f p2+\f{2a+1}2\Big(1+p\big(2q_p(2)+B_{p-1}-B_{p-1}\big(\f
12-a\big)\big)\Big)\mod {p^2},
\endalign$$ where $t'=t$ or $t+1$ according as $\ap<\f p2$ or $\ap>\f p2$. As
$$\f
1{2k-1}=4\f{(2a(a+1)+1)k-a(a+1)}{4k^2-1}-(2a+1)\f{2a+1}{2k+1},$$
from the above and Theorem 4.1 we see that
$$\align &\sum_{k=0}^{(p-1)/2}\f 1{2k-1}\b ak\b {-1-a}k
\\&=4\sum_{k=0}^{(p-1)/2}\b ak\b {-1-a}k\f{(2a(a+1)+1)k-a(a+1)}{4k^2-1}
\\&\q-(2a+1)\sum_{k=0}^{(p-1)/2}\b ak\b {-1-a}k\f{2a+1}{2k+1}
\\&\e 2pt'+2(2a+1)t'\Big(1+p\big(2q_p(2)+B_{p-1}-B_{p-1}\big(\f
12-a\big)\big)\Big)
\\&\q-(-1)^{t'-t}(2a+1)-2(2a+1)t'\big(1+p\big(2q_p(2)-B_{p-1}(-a)+B_{p-1}\big)\big)
\\&=-(-1)^{t'-t}(2a+1)+2pt'+2(2a+1)pt'\Big(B_{p-1}(-a)-B_{p-1}\big(\f
12-a\big)\Big)\mod{p^2}.
\endalign$$
By (3.7),
$$E_{p-2}(-2a)=\f{2^{p-1}}{p-1}\Big(B_{p-1}\big(-a+\f
12\big)-B_{p-1}(-a)\Big)\e B_{p-1}(-a)-B_{p-1}\Big(\f
12-a\Big)\mod{p}.\tag 4.3$$ Now combining the above gives the
result. \pro{Theorem 4.3} Let $p$ be a prime greater than $3$. Then
$$\align
&\sum_{k=0}^{(p-1)/2}\f{\b{2k}k\b{4k}{2k}}{(2k-1)64^k} \e
(-1)^{\f{p+1}2}\f{p+1}2\mod{p^2},\tag i
\\&\sum_{k=0}^{(p-1)/2}\f{\b{2k}k\b{3k}{k}}{(2k-1)27^k}
\e \f 19\Ls p3\big(2^{p+1}-7-6p\big)\mod{p^2},\tag ii
\\&\sum_{k=0}^{(p-1)/2}\f{\b{3k}k\b{6k}{3k}}{(2k-1)432^k}
\e -\f 19\Ls p3\big(2^{p+1}+2+3p\big)\mod {p^2}.\tag iii\endalign$$
\endpro
Proof. It is well known that ([MOS])
$$B_n(1-x)=(-1)^nB_n(x)\qtq{and}E_n(1-x)=(-1)^nE_n(x).$$
For $m=3,4,6$ we see that $\langle -\f 1m\rangle_p=\f{p-1}m$ or
$\f{(m-1)p-1}m$ according as $p\e 1\mod m$ or not.
  Using (4.3) we see that $E_{p-2}\sls 12\e B_{p-1}\sls
14-B_{p-1}\sls 34\e 0\mod p$. Now taking $a=-\f 14$ in Theorem 4.2
and then applying (1.1) yields (i).
 Using (3.10) we see that $E_{p-2}\sls 23=-E_{p-2}\sls 13\e -2q_p(2)\mod
p$. Taking $a=-\f 13$ in Theorem 4.2 and then applying (1.1) yields
(ii). Taking $a=-\f 16$ in Theorem 4.2  and then applying (3.10) and
(1.1) yields (iii). The proof is now complete.

\pro{Theorem 4.4} Let $p>3$ be a prime, $a\in\Bbb Z_p$, $a\not\e
0\mod p$ and $t=(a-\ap)/p$.  Then
$$\sum_{k=0}^{(p-1)/2}\f{\b ak\b{-1-a}k}{2k-1}
\e\cases 2pt-(2a+1)\sum_{k=0}^{p-1}\b{2a}k(-2)^k\mod{p^2} \q\t{if
$\ap<\f p2$,}
\\2p(t+1)-\f{2a+1}{2t+1}
-(2a+1)\big(1+\f 1{2t+1}\big)\sum_{k=0}^{p-1}\b{2a}k(-2)^k\mod{p^2}
\\\qq\qq\qq\qq\qq\q\t{if $\ap>\f p2$ and $t\not\e -\f 12\mod p$,}
\\2a+1+p(1+(2a+1)E_{p-2}(-2a))\mod{p^2}\\\qq\qq\qq\qq\qq\q\t
{if $\ap>\f p2$ and $t\e -\f 12\mod p$.}
\endcases$$
\endpro
Proof. By (1.9),
$$\align \sum_{k=0}^{p-1}\b{2a}k(-2)^k&\e (-1)^{\langle
2a\rangle_p}-(2a-\langle 2a\rangle_p)E_{p-2}(-2a)
\\&=\cases 1-(2a-2\ap)E_{p-2}(-2a)\mod{p^2}&\t{if $\ap<\f p2$,}
\\-1-(2a-2\ap+p)E_{p-2}(-2a)\mod{p^2}&\t{if $\ap>\f p2$.}
\endcases\endalign$$
Thus, from Theorem 4.2 we deduce that
$$\align &\sum_{k=0}^{(p-1)/2}\f{\b ak\b{-1-a}k}{2k-1}+
(2a+1)\sum_{k=0}^{p-1}\b{2a}k(-2)^k
\\&\e\cases 2(a-\ap)\mod{p^2}&\t{if $\ap<\f p2$,}
\\2(p+a-\ap)+(2a+1)pE_{p-2}(-2a)\mod{p^2}&\t{if $\ap>\f p2$.}
\endcases
\endalign$$
This shows that the result is true for $\ap<\f p2$. Now assume
$\ap>\f p2$. From the above we have
$$1+\sum_{k=0}^{p-1}\b{2a}k(-2)^k\e -(2t+1)pE_{p-2}(-2a)\mod
{p^2}.$$ Hence for $t\not\e -\f 12\mod p$,
$$\align &\sum_{k=0}^{(p-1)/2}\f{\b ak\b{-1-a}k}{2k-1}+
(2a+1)\sum_{k=0}^{p-1}\b{2a}k(-2)^k
\\&\e
2(p+a-\ap)-\f{2a+1}{2t+1}\Big(1+
\sum_{k=0}^{p-1}\b{2a}k(-2)^k\Big)\mod{p^2}.
\endalign$$
This yields the result in this case.  For $t\e -\f 12\mod p$, from
the above we have $a-\ap=pt\e -\f p2\mod{p^2}$ and
$\sum_{k=0}^{p-1}\b{2a}k(-2)^k\e -1\mod{p^2}$. Hence
$$\align \sum_{k=0}^{(p-1)/2}\f{\b ak\b{-1-a}k}{2k-1}
&\e 2a+1+2(p+a-\ap)+(2a+1)pE_{p-2}(-2a) \\&\e
2a+1+p(1+(2a+1)E_{p-2}(-2a))\mod{p^2}.\endalign$$ This yields the
result in this case. The proof is now complete.

\par\q
\newline{\bf Acknowledgment}
\par\q
\newline The author is supported by the National
Natural Science Foundation of China (grant No. 11371163).

\end{document}